
\def\today{\ifcase\month\or January\or February\or March\or
April\or May\or June\or July\or August\or September\or
October\or November\or December\fi \space\number\day,
\number\year}



\font \bbrm=cmbx10 at 12pt

\font \ninerm= cmr10 at 9pt
\def\smalltype{\ninerm}

\def\bigtype{\bbrm}

\hsize=13.5cm
\magnification=1200
\def\ce{\centerline}

\def\title #1{\null\bigskip\ce{\bigtype #1}\bigskip}






\font\tenboldgreek=cmmib10  \font\sevenboldgreek=cmmib10 at
7pt
\font\fiveboldgreek=cmmib10 at 7pt
\newfam\bgfam
\textfont\bgfam=\tenboldgreek \scriptfont\bgfam=\sevenboldgreek
\scriptscriptfont\bgfam=\fiveboldgreek


\font\tengerman=eufm10 \font\sevengerman=eufm7 \font\fivegerman=eufm5
\font\tendouble=msym10 \font\sevendouble=msym7 \font\fivedouble=msym5

\textfont4=\tengerman \scriptfont4=\sevengerman
\scriptscriptfont4=\fivegerman
\newfam\dbfam
\textfont\dbfam=\tendouble \scriptfont\dbfam=\sevendouble
\scriptscriptfont\dbfam=\fivedouble

\mathchardef\ng="702D
\mathchardef\dbA="7041
\mathchardef\sm="7072
\mathchardef\nvdash="7030
\mathchardef\nldash="7031
\mathchardef\lne="7008
\mathchardef\sneq="7024
\mathchardef\spneq="7025
\mathchardef\sne="7028
\mathchardef\spne="7029
\mathchardef\ltms="706E
\mathchardef\tmsl="706F

\mathchardef\dbA="7041


\mathchardef\dbA="7041 
\mathchardef\dbB="7042 
\mathchardef\dbC="7043 
\mathchardef\dbD="7044 
\mathchardef\dbE="7045 
\mathchardef\dbF="7046 
\mathchardef\dbG="7047 
\mathchardef\dbH="7048 
\mathchardef\dbI="7049 
\mathchardef\dbJ="704A 
\mathchardef\dbK="704B 
\mathchardef\dbL="704C 
\mathchardef\dbM="704D 
\mathchardef\dbN="704E 
\mathchardef\dbO="704F 
\mathchardef\dbP="7050 
\mathchardef\dbQ="7051 
\mathchardef\dbR="7052 
\mathchardef\dbS="7053 
\mathchardef\dbT="7054 
\mathchardef\dbU="7055 
\mathchardef\dbV="7056 
\mathchardef\dbW="7057 
\mathchardef\dbX="7058 
\mathchardef\dbY="7059 
\mathchardef\dbZ="705A

\def\sdp{\times \hskip -0.3em {\raise 0.3ex
\hbox{$\scriptscriptstyle |$}}} 


\def\const{\mathop{\rm const}}

\def\grad{\mathop {\rm grad}\nolimits}




\def\ot{{\overline t}}

\def\oV{{\overline V}}

\def\ox{{\overline x}}

\def\oy{{\overline y}}








\def\sqr#1#2{{\vcenter{\hrule height.#2pt\hbox{\vrule
width.#2pt height#1pt \kern#1pt \vrule width.#2pt}\hrule
height.#2pt}}}
\def\square{\mathchoice{\sqr34}{\sqr34}{\sqr{2.1}3}{\sqr{1.5}3}}

\def\buildrul#1\under#2{\mathrel{\mathop{\null#2}\limits_{#1}}}

\def\boxit#1{\vbox{\hrule\hbox{\vrule\kern3pt\vbox{\kern3pt#1
\kern3pt}\kern3pt\vrule}\hrule}}

\def\hb{\hfill\break}
\def\subheading#1{\medskip\goodbreak\noindent{\bf
#1.}\quad}
\def\sect#1{\goodbreak\bigskip\centerline{\bf#1}\medskip}
\def\pr{\smallskip\noindent{\bf Proof:\quad}}

\def\longmapright #1 #2 {\smash{\mathop{\hbox to
#1pt {\rightarrowfill}}\limits^{#2}}}
\def\longmapleft #1 #2 {\smash{\mathop{\hbox to
#1pt {\leftarrowfill}}\limits^{#2}}}

\def\ref#1{\par\smallskip\hang\indent\llap{\hbox
to \parindent{#1\hfil\enspace}}\ignorespaces}

\def\back{{\raise 2.5pt\hbox{$\,\scriptscriptstyle\backslash\,$}}}
\def\part{\partial}
\def\lwr #1{\lower 5pt\hbox{$#1$}\hskip -3pt}
\def\rse #1{\hskip -3pt\raise 5pt\hbox{$#1$}}
\def\lwrs #1{\lower 4pt\hbox{$\scriptstyle #1$}\hskip -2pt}
\def\rses #1{\hskip -2pt\raise 3pt\hbox{$\scriptstyle #1$}}

\def\<#1{\left\langle{#1}\right\rangle}

\def\subinbn{{\subset\hskip-8pt\raise
0.95pt\hbox{$\scriptscriptstyle\subset$}}}

\def\llvdash{\mathop{\|\hskip-2pt \raise 3pt\hbox{\vrule
height 0.25pt width 1.5cm}}}

\def\lvdash{\mathop{|\hskip-2pt \raise 3pt\hbox{\vrule
height 0.25pt width 1.5cm}}}

\def\fakebold#1{\leavevmode\setbox0=\hbox{#1}%
  \kern-.025em\copy0 \kern-\wd0
  \kern .025em\copy0 \kern-\wd0
  \kern-.025em\raise.0333em\box0 }

\font\msxmten=msxm10
\font\msxmseven=msxm7
\font\msxmfive=msxm5
\newfam\myfam
\textfont\myfam=\msxmten
\scriptfont\myfam=\msxmseven
\scriptscriptfont\myfam=\msxmfive
\mathchardef\ggarrow="7010
\mathchardef\rhookupone="7016
\mathchardef\BECAUSE="702A
\mathchardef\ldh="700D
\mathchardef\leg="7053
\mathchardef\ANG="705E
\mathchardef\lcu="7070
\mathchardef\rcu="7071
\mathchardef\leseq="7035
\mathchardef\qeeg="703D
\mathchardef\qeel="7036
\mathchardef\blackbox="7004
\mathchardef\bbx="7003
\mathchardef\simsucc="7025

\def\bigsquare{{\fam=\myfam\bbx}}

\font\tencaps=cmcsc10
\def\smallcaps{\tencaps}

\def\author#1{\bigskip\ce{\smallcaps #1}\medskip}

\def\dotsq{\ooalign{\hfil\raise
1pt\hbox{$\cdot$}\hfil\cr\cr\hbox{$\bigsquare$}}}

\def\frac #1#2{{#1\over #2}}

\def\endproclaim{\par\rm\ifdim\lastskip<\medskipamount\removelastskip
\penalty55\medskip\fi}

\def\leaderfill{\leaders\hbox to 1em{\hss.\hss}\hfill}

\def\upddots{\mathinner{\mkern
1mu\raise 1pt \hbox{.}\mkern 2mu \mkern
2mu \raise 4pt\hbox{.}\mkern 1mu \raise 7pt\vbox {\kern 7
pt\hbox{.}}} }

\def\sqr#1#2{{\vcenter{\hrule height.#2pt\hbox{\vrule
width.#2pt height#1pt \kern#1pt \vrule width.#2pt}\hrule
height.#2pt}}}
\def\square{\mathchoice{\sqr34}{\sqr34}{\sqr{2.1}3}{\sqr{1.5}3}}

\def\buildrul#1\under#2{\mathrel{\mathop{\null#2}\limits_{#1}}}

\def\boxit#1{\vbox{\hrule\hbox{\vrule\kern3pt\vbox{\kern3pt#1
\kern3pt}\kern3pt\vrule}\hrule}}

\def\subheading#1{\medskip\goodbreak\noindent{\bf
#1.}\quad}
\def\sect#1{\goodbreak\bigskip\centerline{\bf#1}\medskip}
\def\pr{\smallskip\noindent{\bf Proof:\quad}}

\def\longmapright #1 #2 {\smash{\mathop{\hbox to
#1pt {\rightarrowfill}}\limits^{#2}}}
\def\longmapleft #1 #2 {\smash{\mathop{\hbox to
#1pt {\leftarrowfill}}\limits^{#2}}}

\def\ref#1{\par\smallskip\hang\indent\llap{\hbox
to \parindent{#1\hfil\enspace}}\ignorespaces}

\def\back{{\raise 2.5pt\hbox{$\,\scriptscriptstyle\backslash\,$}}}
\def\part{\partial}
\def\lwr #1{\lower 5pt\hbox{$#1$}\hskip -3pt}
\def\rse #1{\hskip -3pt\raise 5pt\hbox{$#1$}}
\def\lwrs #1{\lower 4pt\hbox{$\scriptstyle #1$}\hskip -2pt}
\def\rses #1{\hskip -2pt\raise 3pt\hbox{$\scriptstyle #1$}}

\def\<#1{\left\langle{#1}\right\rangle}

\def\subinbn{{\subset\hskip-8pt\raise
0.95pt\hbox{$\scriptscriptstyle\subset$}}}

\def\llvdash{\mathop{\|\hskip-2pt \raise 3pt\hbox{\vrule
height 0.25pt width 1.5cm}}}

\def\lvdash{\mathop{|\hskip-2pt \raise 3pt\hbox{\vrule
height 0.25pt width 1.5cm}}}

\def\fakebold#1{\leavevmode\setbox0=\hbox{#1}%
  \kern-.025em\copy0 \kern-\wd0
  \kern .025em\copy0 \kern-\wd0
  \kern-.025em\raise.0333em\box0 }

\font\msxmten=msxm10
\font\msxmseven=msxm7
\font\msxmfive=msxm5
\newfam\myfam
\textfont\myfam=\msxmten
\scriptfont\myfam=\msxmseven
\scriptscriptfont\myfam=\msxmfive
\mathchardef\rhookupone="7016
\mathchardef\ldh="700D
\mathchardef\leg="7053
\mathchardef\ANG="705E
\mathchardef\lcu="7070
\mathchardef\rcu="7071
\mathchardef\leseq="7035
\mathchardef\qeeg="703D
\mathchardef\qeel="7036
\mathchardef\blackbox="7004
\mathchardef\bbx="7003
\mathchardef\simsucc="7025

\def\bigsquare{{\fam=\myfam\bbx}}

\font\tencaps=cmcsc10
\def\smallcaps{\tencaps}

\def\author#1{\bigskip\ce{\smallcaps #1}\medskip}

\def\frac #1#2{{#1\over #2}}

\def\upddots{\mathinner{\mkern
1mu\raise 1pt \hbox{.}\mkern 2mu \mkern
2mu \raise 4pt\hbox{.}\mkern 1mu \raise 7pt\vbox {\kern 7
pt\hbox{.}}} }
\def\circle{{\raise 1pt\hbox{${\scriptstyle\bigcirc}$}}}


\def\ce{\centerline}
\def\sect#1{\goodbreak\bigskip\centerline{\bf#1}\medskip}

\def\hb{\hfill\break}

\nopagenumbers

\vsize=24.5truecm
\null
\sect{ A NEW APPROACH TO INVESTIGATION OF EVOLUTION}
\vskip-0.5truecm
\sect{DIFFERENTIAL EQUATIONS IN BANACH SPACES}
\sect{Ya.I. Alber\footnote*{\smalltype This
research is supported in part by the Ministry of Scence Grant
3481-1-91 and by the Ministry of Absorption, Center for Absorption in
Science.}}
\centerline{Department of Mathematics }
\centerline{Technion -- Israel Institute of Technology }
\ce{Haifa 32000, Israel}
\vskip0.75truecm
\sect{Introduction}
\ospace

Known investigations of nonlinear evolution equations
$$
{dx\over dt} + A(t)x(t) = f(t)\ ,\quad x(t_{0}) = x^{0}\
,\quad t_{0} \le  t < \infty\ ,
\eqno(0.1)$$
with monotone operators $A(t)$ acting from reflexive Banach space $B$
to  dual
space $B^*$, usually assume that along with $B$ and $B^*$ there is
a Hilbert  space
$H$ and continuous imbedding $B \hookrightarrow  H$ in the triplet
$$
B \hookrightarrow  H \hookrightarrow  B^*\ ;
\eqno(0.2)$$
and that $B$ is dense in $H$.
 The  stabilization  of  solutions  of  evolution
equations has been proven either in the sense of weak
convergence in $B$  or
in the norm of $H$ space, and only asymptotic estimates of stabilization
rate
have been obtained [15].
In the present paper  we  consider  equations  of
type (0.1) without conditions (0.2) and establish stabilization
with  both
asymptotic and nonasymptotic rate estimated in the norm of $B$.
Our research
is based on a new technique  using  Banach  space  geometry,
parallelogram
inequalities,  estimates  of   duality   mappings,   nonstandard
Lyapunov
functionals in Banach space, and estimates  of  solutions  to
differential
inequalities.
\bigskip

\sect{1. Differential Inequalities}

In this section differential inequalities
$$
{d\lambda (t)\over dt} \le  - \alpha (t)\psi \big(\lambda (t)\big)
+ \gamma (t)\ ,\quad \lambda (t_{0}) = \lambda _{0}\ ,\quad t \ge
 t_{0}
\eqno(1.1)$$
are investigated for a nonnegative function $\lambda (t)$.
 Such  inequalities  are
the most powerful tool, unique in a lot of cases, to prove convergence  and
stability and to establish the estimates of convergence rate
solutions  of evolution equations.
The results that will be presented play an  auxiliary
role. However, they have their own importance as well (see also [16]).

\headline{\hfill\smalltype\folio\hfill}
\proclaim Lemma 1.
If,  in  inequality $(1.1)$, $\psi (\lambda )$   is  a  positive
 continuous
function for all $\lambda  > 0$, $\psi (0) = 0$ and, in addition,
$$
\lim_{t\rightarrow \infty } {\gamma (t)\over \alpha (t)} = 0\ ,
$$
$$
\int^{\infty }_{t_{0}}\alpha (\tau)d\tau = \infty \ ,
$$
then $\lambda (t) \rightarrow  0$ as $t \rightarrow  \infty$.

\pr
Let us consider the alternative
$$
H_{1}: \psi \big(\lambda (t)\big) < S(t)\ ,\quad H_{2}: \psi \big(
\lambda (t)\big) \ge  S(t)\ ,
$$
$$
S(t) = {1\over \int^{t}_{t_{0}}\alpha (\tau)d\tau} + {\gamma
(t)\over \alpha (t)}\ .
$$
Introduce the sets
$$
\eqalignno{&T^{i}_{1} =\big\{t_{0} \le  t \in  [t_{i},\ot_{i}]
\subseteq  R^{+}: H_{1}\hbox{ is true }\big\}\ ,\quad T_{1} =
\bigcup T^{i}_{1}\ ,\quad i = 1,2,\ldots &(1.2)\cr
&T^{j}_{2} =\big\{t_{0} \le  t \in  [t_{j},\ot_{j}\big\}] \subseteq
 R^{+}: H_{2}\hbox{ is true }\big\}\ ,\quad T_{2} = \bigcup
T^{j}_{2}\ ,\quad j = 1,2,\ldots
 .  &(1.3)\cr}$$
It is obvious that $T^{i}_{1}$ and $T^{j}_{2}$ alternate and
$T = [t_{0},\infty ) = T_{1} \cup  T_{2}$.   We  will
show $T_{1}$ to be unbounded.  Let the opposite be true.
Then  there  is  some
$t=\tau_{1}$ such that  for  all $t \ge  \tau_{1}$  the
hypothesis $H_{2}$  occurs  and  the
differential inequality
$$
{d\lambda (t)\over dt} \le  - {\alpha (t)\over
\int^{t}_{t_{0}}\alpha (\tau)d\tau}
\eqno(1.4)$$
is realized.  Hence,
$$
\lambda (t) \le  \lambda (\tau_{1}) -
\int^{t}_{\tau_{1}}{\alpha (\tau)\over A(\tau)}d\tau\ ,
\quad A(\tau) = \int^{\tau}_{t_{0}}\alpha (s)ds\ .
$$
By virtue of the Cauchy integral criterion,
$$
\int^{\infty }_{\tau_{1}}{\alpha (\tau)\over A(\tau)}d\tau = \infty
$$
because
$$
\int  {\alpha (t)\over A(t)}d\tau = \ell nA(t) \rightarrow
\infty\quad \hbox{ as }\quad t \rightarrow  \infty \ .
$$
It follows from this that, beginning with some $t = \tau_{2}$,
the  function $\lambda (t)$
becomes negative, but this is in contradiction with the conditions
of  the
lemma.
Thus $\lambda (t)$  is  estimated  from  above  by  a  function
 which  is
monotonically decreasing and vanishes  on  the  set $T_{1}$.
Because  of  the
$H_{2}$-hypothesis, $\lambda (t)$ is also estimated on every
interval $T^{j}_{2}$ by a  decreasing
function of the (1.4) kind.  The lemma is proved.

\subheading{Remark 1}
  The assertion of Lemma 1 for the linear case
$\big(\psi (\lambda )=\lambda \big)$ has
been obtained previously in [1], using the formula
$$
\lambda (t) \le  \lambda (t_{0})e^{-\int^{t}_{t_0}\alpha (\tau)d\tau}
+ \int^{t}_{t_{0}}\gamma (\tau)e^{-\int^{t}_{\tau}\alpha (s)ds}d\tau\ .
$$
Now we will introduce stronger requirements concerning functions
$\psi (\lambda ),
\alpha (t)$ and $\gamma (t)$, that enable us to prove stabilization
of $\lambda (t)$  to  zero,  as
well as  to  establish  both  asymptotic  and  nonasymptotic  estimates  of
the stabilization rate.

First of all we will formulate two general statements about
inequality
(1.1).
Let $\gamma (t)$ and $\alpha (t)$  be  continuous,  nonnegative
and  nonincreasing
functions, $\psi (\lambda )$ a continuous positive  strictly
increasing  function  with
$\psi (0) = 0$, and let $F(t)$ and $\phi (\lambda )$ be
antiderivatives of  the  functions $\alpha (t)$
and $1/\psi (\lambda )$, respectively.
These  functions  are  assumed  to  be  defined,
$\phi ^{-1}(z)$ exists and is single-valued on the corresponding
set $Z$, and $F(t) > 0$
for all $t \ge  t_{0}$.

We introduce the following notations [4]:
\item{(1)}$v(t) = \psi ^{-1}\big(c_{0}{\gamma (t)\over \alpha
(t)}\big)$,
$t \ge  t_{0}$, where $c_{0} > 1$ is  an  arbitrary  parameter
and $\psi ^{-1}(\cdot )$ is the function inverse to $\psi (\cdot )$,
$v_{0} = v(t_{0})$;

$$
 u(t,C) = \phi ^{-1}\big[\phi (C)-a(F(t)-F(t_{0}))\big]\ ,\quad
 a = c^{-1}_{0}(c_{0}-1) > 0\ ,\quad C \ge  0\ .
\leqno(2)$$

Next, we reduce the main property $(P)$:  on  the  interval
$[t_{0},\infty )$  the
functions $v(t)$ and
$$
w\big(t,v(z)\big) = \phi ^{-1}\big[\phi (v(z))-a(F(t)-F(z))\big]\ ,
$$
where $z \in  [t_{0},\infty )$ is fixed, either coincide completely
or  intersect  at  no
more  than  two  points  (the  points  of  tangency  are  not
 regarded  as
intersection points of the functions except, possibly, for $t
= t_{0})$.

\proclaim Lemma 2.
Suppose that $(1)$ $(P)$ holds for $z = t_{0}$, $(2)$ $u(t,v_{0})
> v(t)$  as
$t\rightarrow \infty $, and $(3)$ $c_{0}$ is chosen
such that $u(t,v_{0}) \ge  v(t)$ as $t \rightarrow  t_{0}$.
Then  the
solution of inequality $(1.1)$ $\lambda (t)\rightarrow 0$ and for
all $t \ge  t_{0}$ satisfies  the
next estimate
$$
\lambda (t) \le  u(t,C)\ ,\quad C = \max \{\lambda _{0},v_{0}\}\  .
\eqno(1.5)$$

\subheading{Remark 2}
Condition 2 is understood in the sense that  there  exists
the number $N$ such that $u(t,v_{0}) \ge  v(t)$ for all $t \ge  N$.
This also applies  to
similar conditions contained in Lemmas 2 and 3.

\pr
 Consider the following alternative
$$
\eqalignno{&H_{1}: \lambda (t) < \psi ^{-1}\big(c_{0}\gamma
(t)/\alpha (t)\big)
&(1.6)\cr
&H_{2}: \lambda (t) \ge  \psi ^{-1}\big(c_{0}\gamma (t)/\alpha (t)
\big)\cr}
$$
and the intervals $T^{i}_{1}$ and $T^{j}_{2}$ with respect to
(1.2) and (1.3).

Assume, at first, that $T_{1} = \phi $.
In that case, for all $t \in  T = [t_{0},\infty )$,
we have $\gamma (t) \le  c^{-1}_{0}\alpha (t)\psi (\lambda (t))$
and inequality
$$
{d\lambda (t)\over dt} \le  -a\alpha (t)\psi (\lambda (t))\ ,
\quad \lambda (t_{0}) = \lambda _{0}\ .
$$
Its solution is estimated by the known formula
$$
\lambda (t) \le  \phi ^{-1}\big[\phi
(\lambda _{0})-a(F(t)-F(t_{0}))\big]\ ,\quad t \ge  t_{0}\  .
\eqno(1.7)$$
Now let $T_{2}= \emptyset$.  Then, estimate (1.6) is true.

If $T_{1}$ or $T_{2}$ is  finite,  then  the  estimates  (1.6)
or  (1.7)  hold,
starting from some $\widetilde{t} \ge  t_{0}$, and a situation
arising under $t \le  \widetilde{t}$ is described
below.

Finally, assume $T_{1}$ and $T_{2}$ to be infinite sets at the
same  time.   Let
us choose an arbitrary interval $T^{i}_{1} \in  T_{1}$.
This determines the next interval
$T^{i}_{2}$ (if $t_{0} \in  T_{1}$) or $T^{i+1}_{2}$ (if $t_{0}
\in  T_{2}$).
It  is  obvious  that  (1.6)  is
fulfilled for all $t \in  T^{i}_{1}$, but that on the interval
$T^{i}_{2}$
$$
\lambda (t) \le  \phi ^{-1}\bigg[\phi\big(\lambda (\ot_{i})\big
)-a \int^{t}_{\ot_{i}}\alpha (s)ds\bigg]\ .
\eqno(1.8)$$
One should consider the following two cases:
\item{(a)} $\lambda _{0} \ge  v(t_{0})$.
Here $t_{0} \in  T_{2}$ and the function $u(t,\lambda _{0})$
is  majorant.
Indeed, the function $u(t,\lambda _{0})$ cannot be intersected
with $v(t)$ at any value $t$
because, in the opposite case, at least three points of
intersection  that
do not match the property $(P)$ would appear.
\item{(b)}$\lambda _{0} < v(t_{0})$.
Let  us  construct $u\big(t,v(t_{0})\big)$.  According  to  the
assumptions of the lemma, $u(t,v_{0})$ majorizes $v(t)$ as
$t_{0}\leftarrow  t \rightarrow  \infty $  and  by
virtue of property
$$
u(t,R_{1}) \ge  u(t,R_{2})\ ,\quad \forall\,  R_{1} \ge  R_{2}
\eqno(1.9)$$
it  majorizes $u(t,\lambda _{0})$  along  with  function
$\phi ^{-1}(\cdot )$  in  (1.8),   because
$\lambda (\ot_{i})=v(\ot_{i})$.
We remark that $u(t,\lambda _{0})$ has to intersect $v(t)$ only
once,  so
$i = 1$.
Hence, the general estimate is described by the formula (1.5).  The
lemma is proven.

\proclaim Lemma 3.
Suppose that $(1)$ property $(P)$ holds  for  all $z \in
 [t_{0},\infty )$;
$(2)$ $u(t,v_{0}) \le  v(t)$ as $t \rightarrow  \infty $.
Then $\lambda (t) \rightarrow  0$ and  the  following  assertions
are true:
\item{\rm (i)} If $\lambda _{0} \le  v_{0}$ and $c_{0}$ is chosen
so that $u(t,v_{0}) \le  v(t)$ as $t \rightarrow  t_{0}$,  then,
for the solution of inequality $(1.1)$ the estimate
$$
\lambda (t) \le  v(t)\ ,\quad t \ge  t_{0}
$$
is valid;
\item{\rm  (ii)} In all the remaining cases,
$$
\eqalignno{&\lambda (t) \le  u(t,C)\ ;\quad C = \max
\big\{\lambda _{0},v(t_{0})\big\}\ ,\quad t_{0} \le  t \le \ot&
(1.10)\cr
& \lambda (t) \le  v(t)\ ,\quad t \ge \ot&(1.11)\cr}$$
where $\ot$  is the unique root of  the  equation $u(t,C) = v(t)$
on  the
interval $[t_{0},\infty )$.

\pr
  Let the assumption (i) be valid.  By (2) and property $(P)$ the
function $u(t,v_{0})$ is no greater than $v(t)$ for all $t \ge
 t_{0}$,  otherwise  they
would have no less than three intersection points on the
interval $[t_{0},\infty )$.
If $\lambda _{0} \le  v_{0}$, then the function $w(z,v(z))$ is
majorized also by  the  function
$v(t)$, since, in the opposite case, either three  intersection
 points  must
appear, or one intersection point $t > z$  and  one  tangent
point $t  =  z$.
Neither is appropriate by virtue of property (1.9) so long as
 the  tangent
point is disintegrated into two intersection points by a tiny
perturbation.
So function $v(t)$ is a majorant.

Now let $\lambda _{0} \le  v_{0}$ but $u(t,v_{0}) > v(t)$  as
$t \rightarrow  t_{0}$.  It is  clear  that  in
that variant exactly one intersection point for $v(t)$  and
$v(t,v_{0})$  exists.
The same applies to $\lambda _{0} > v_{0}$.
The  major  estimates  are  determined  by
relations (1.10) and (1.11).  The lemma is proven.

We shall reduce  the  most  important  versions  of  Lemmas  2
 and  3
(cf.~[10]).

\proclaim Lemma 4.
 In inequality $(1.1)$ let $\lambda (t) \ge  0$, $\psi (\lambda )
= \lambda $,
$$
\alpha (t) = {b\over t}\ ,\quad \gamma (t) = {d\over t^{n}}\ ,
\quad n > 1\ ,\quad b > 0\ ,\quad d \ge  0\ .
$$
Then $\lambda (t) \rightarrow  0$ as $t \rightarrow  \infty $ and,
also:
\item{$ (1)$} Assume $ab > n-1$, i.e., $c_{0} > b/(b-n+1)$.
Then, the following are true:
\itemitem{\rm (i)}If $\lambda _{0} \le  c_{0}db^{-1}t^{1-n}_{0}$,
then, for all $t \ge  t_{0}$,
$$
\lambda (t) \le  {c_{0}d\over b}\left({1\over t}\right)^{n-1}\ ;
$$
\itemitem{\rm (ii)} If $\lambda _{0} \ge  c_{0}db^{-1}t^{1-n}_{0}$,
then, the estimates
$$
\eqalign{&\lambda (t) \le  \lambda_0\left({t_{0}\over t}\right)
^{ab}\ ,\quad t_{0} \le  t \le \ot\cr
&\lambda (t) \le  {c_{0}d\over b}\left({1\over t}\right)^{n-1}\ ,
\quad t >\ot\ ,\quad \ot= (c^{-1}_{0}d^{-1}b\lambda _{0}t^{ab}_{0}
)^{1/(ab+1-n)}\cr}
$$
are satisfied.
\item{$(2)$} Assume $ab \le  n-1$, i.e., $c_{0} \le  b/(b-n+1)$.
Then, for all $t \ge  t_{0}$
$$\lambda (t) \le  C\left({t_{0}\over t}\right)^{ab}\ ,\quad
 C = \max \left\{\lambda _{0}, {c_{0}d\over b}
t^{1-n}_{0}\right\}\ .$$

\proclaim Lemma 5.
In inequality $(1.1)$ let $\lambda (t) \ge  0$, $\psi (\lambda )
= \lambda $,
$$
\alpha (t) = {b\over t^m}\ ,\quad 0 \le  m < 1\ ,\quad \gamma
(t) = {d\over t^n}\ ,\quad n > m\ ,\quad b > 0\ ,\quad d \ge  0\ .
$$
 Then $\lambda (t) \rightarrow  0$ as $t \rightarrow  \infty $ and,
also:
\item{$(1)$} If $\lambda _{0} \le  c_{0}db^{-1}t^{m-n}_{0}$ and
 $ab \ge  t^{m-1}_{0}(n-m)$, i.e., $c_{0} \ge
 bt^{1-m}_{0}/(bt^{1-m}_{0}-n+m)$,
then, for all $t \ge  t_{0}$,
$$
\lambda (t) \le  {c_{0}d\over b}\left({1\over t}\right)^{n-m}\ ;
$$
\item{$(2)$} In all the remaining cases
$$
\eqalign{&\lambda (t) \le  C~\exp \left[-{ab\over
1-m}(t^{1-m}-t^{1-m}_{0})\right]\ ,\quad t_{0} \le  t \le \ot
\cr
&
\lambda (t) \le  {c_{0}d\over b}\left({1\over t}\right)^{n-m}\ ,
\  t >\ot\ ,\  C = \max \left\{\lambda _{0},{c_{0}d\over
b} t^{m-n}_{0}\right\}\ ,\qquad\cr}
$$
where $\ot$  is the unique root of equation
$$
C \exp \left[{ab\over m-1}(t^{1-m}-t^{1-m}_{0})\right] =
{c_{0}d\over b}\left({1\over t}\right)^{n-m}\  .
$$

The asymptotic estimate is
$$
\lambda (t) = O(t^{-n})
$$
as $m = 0.$

\proclaim Lemma 6.
In inequality $(1.1)$ let $\lambda  \ge  0$, $\psi (\lambda )
= \lambda ^{\nu }$, $\nu  > 1,$
$$
\alpha (t) = {b\over t}\ ,\quad \gamma (t) = {d\over t^{n}}\ ,
\quad n > 1\ ,\quad b > 0\ ,\quad d \ge  0\ .
$$
Then $\lambda (t) \rightarrow  0$  as $t \rightarrow  \infty $
and, if
$$
ab\left({dc_{0}\over b}\right)^{\nu-1\over \nu }<{n-1\over\nu}
 t^{{(n-1)(\nu -1)\over \nu }}_{0}
\eqno(1.12)$$
holds, then, for all $t \ge  t_{0}$,
$$
\lambda (t) \le  \big[C^{1-\nu }+(\nu -1)ab\ln
t/t_{0}\big]^{-1/(\nu -1)}
$$
where $C = \max \big\{\lambda _{0},(c_{0}db^{-1})^{1/\nu
}t^{(n-1)/\nu }_{0}\big\}$.

\proclaim Lemma 7.
In inequality $(1.1)$ let $\lambda (t) \ge  0$, $\psi (\lambda
) = \lambda ^{\nu }$, $\nu  > 1,$
$$\alpha (t) = {b\over t^{m}}\ ,\quad 0 \le  m < 1\ ,\quad
\gamma (t) = {d\over t^{n}}\ ,\quad n > m\ ,\quad b > 0\ ,\quad
 d \ge  0\ .
$$
Then $\lambda (t) \rightarrow  0$ as $t \rightarrow  \infty $ and,
also:
\item{$(1)$}
Assume $\nu ^{-1}(n-m) < (\nu -1)^{-1}(1-m)$ and
\itemitem{\rm (i)}if  $\lambda _{0} \le  (c_{0}d b^{-1})^{1/\nu }
t^{(m-n)/\nu }_{0}$ and
$$
ab\left({c_{0}d\over b}\right)^{(\nu-1)/\nu}\ge{n-m\over\nu}t ^{m-p}_0\
,\quad p = 1 - {(\nu -1)(n-m)\over \nu }
\eqno(1.13)$$
then, for all $t \ge  t_{0}$
$$
\lambda (t) \le  \left({c_{0}d\over b}\right)^{1/\nu }\left({1\over
t}\right)^{(n-m)/ \nu }\ .
$$
\itemitem{\rm (ii)} In all the remaining cases
$$
\eqalignno{&\lambda (t) \le  \left[C^{1-\nu }+ab{\nu -1\over 1-m}(
t^{1-m}-t^{1-m}_{0})\right]^{-1/(\nu -1)}\ ,\quad t_{0} \le  t \le
\ot\cr
&\lambda (t) \le  \left({c_{0}d\over b}\right)^{1/\nu }\left(
{1\over t}\right)^{(n-m)/ \nu }\ ,\ t >\ot\ ,\ C =
\max \left\{\lambda _{0},\left({c_{0}d\over b}\right)^{1/\nu }
t^{(m-n)/ \nu }_{0}\right\}\ ,\qquad\qquad
&(1.14)\cr}$$
where $\ot$  is a unique root of the equation
$$
\left[C^{1-\nu }+ab{\nu -1\over 1-m}(t^{1-m}-t^{1-m}_{0})\right]
^{-1/(\nu -1)} = \left({c_{0}d\over b}\right)^{1/\nu }\left(
{1\over t}\right)^{(n-m)/ \nu }\ .
$$
\item{$(2)$} Assume $(1-m)(\nu -1)^{-1} \le  (n-m)\nu ^{-1}$.
If an inequality opposite  to  $(1.13)$
occurs, then, for all $t > t_{0}$,
$$
\lambda (t) \le  \left[C^{1-\nu }+ab{\nu -1\over 1-m}(t^{1-m}
-t^{1-m}_{0})\right]^{-1/(\nu -1)}
$$
where constant $C$ coincides with $(1.14)$.

The asymptotic estimate is
$$
\lambda (t) = O(t^{-p})\ ,\quad p = \min \left\{{n\over\nu},{1\over \nu
-1}\right\}
$$
as $m = 0.$

\proclaim Lemma 8.
In inequality $(1.1)$ let $\lambda (t) \ge  0$, $\psi (\lambda )
= \lambda ^{\nu }$, $0 < \nu  < 1$,
$$
\alpha (t) = {b\over t}\ ,\quad \gamma (t) = {d\over t^{n}}\ ,\quad
n > 1\ ,\quad b > 0\ ,\quad d \ge  0\ .
$$
Then $\lambda (t) \rightarrow  0$ as $t \rightarrow  \infty $ and,
also:
\item{$(1)$}
If $\lambda _{0} < (c_{0}db^{-1})^{1/\nu } t^{(1-n)/\nu }_{0}$ and
 $c_{0}$ is chosen from $(1.12)$,  then,  for
all $t \ge  t_{0}$,
$$
\lambda (t) \le\left({c_{0}d\over b}\right)^{1/\nu }\left({1\over t}
\right)^{(n-1)/ \nu }\ .
$$
\item{$(2)$} In all the remaining cases the estimates
$$
\eqalign{&\lambda (t) \le  C\left[1-C^{\nu -1}(1-\nu )ab \ln
{t\over t_{0}}\right]^{1/(1-\nu )}\ ,\quad t_{0} \le  t \le \ot\
,\cr
&\lambda (t) \le  \left({c_{0}d\over b}\right)^{1/\nu }\left(
{1\over t}\right)^{(n-1)/ \nu }\ ,\qquad t >\ot \ ,\cr}
$$
are satisfied, where $C = \max \big\{\lambda _{0},(c_{0}db^{-1})
^{1/\nu }t^{(1-n)/\nu }_{0}\big\}$ and $\ot$  is a unique root
of equation
$$
\left({c_{0}d\over b}\right)^{1/\nu }\left({1\over t}\right)^{(
n-1)/ \nu } = C\left[1-(1-\nu )C^{\nu -1}ab \ln{t\over t_{0}}
\right]^{1/(1-\nu )}
$$
 belonging to the interval $\big[t_{0},t_{0}\exp
\big\{C^{1-\nu }/ab(1-\nu )\big\}\big]$.

\proclaim Lemma 9.
In inequality $(1.1)$  let $\lambda (t) \ge  0$, $\psi (\lambda )
= \lambda ^{\nu }$, $0 < \nu  < 1$,
$$
\alpha (t) = {b\over t^{m}}\ ,\quad 0 \le  m < 1\ ,\quad \gamma
(t) = {d\over t^{n}}\ ,\quad n > m\ ,\quad b > 0\ ,\quad d \ge  0\ .
$$
 Then $\lambda (t) \rightarrow  0$  as $t \rightarrow  \infty $
and, also:
\item{$(1)$} If $\lambda _{0} \le  (c_{0}db^{-1})^{1/\nu }
t^{(m-n)/\nu }_{0}$  and $c_{0}$  is chosen from  $(1.13)$,
 then,  for all $t \ge  t_{0}$,
$$
\lambda (t) \le   \left({c_{0}d\over b}\right) ^{1/\nu } \left(
{1\over t}\right)^{(n-m)/ \nu }\ .
$$
\item{$(2)$} In all the remaining cases
$$
\eqalign{&\lambda (t) \le  C\left[1+{1-\nu \over m-1}ab\,
C^{\nu -1}(t^{1-m}-t^{1-m}_{0})\right]^{1/(1-\nu )}\ ,\quad
t_{0} \le  t \le \ot\ ,\cr
&\lambda (t) \le   \left({c_{0}d\over b}\right)^{1/\nu }\left(
{1\over t}\right)^{(n-m)/ \nu }\ ,\quad  t >\ot\ ,\cr}
$$
 where $C = \max \big\{\lambda _{0},(c_{0}db^{-1})^{1/\nu }
t^{(m-n)/\nu }_{0}\big\}$ and  $\ot$ is  a  unique  root
 of  the
equation
$$
\left[C^{\nu -1}+{1-\nu \over m-1}ab (t^{1-m}-t^{1-m}_{0})\right]
^{1/(1-\nu )} =\left({c_{0}d\over b}\right)^{1/\nu}\left({1\over t}
\right)^{(n-m)/ \nu }
$$
 belonging to interval $\big[ t_{0},\big(C^{1-\nu }(1-m)(1-\nu )
^{-1}(ab)^{-1}+t^{1-m}_{0}\big)^{1/(1-m)}\big]   $.

\medskip
The asymptotic estimate is
$$
\lambda (t) = O(t^{-n/\nu })\ ,
$$
as $m = 0.$

\proclaim Lemma 10.   Let $\lambda (t) \ge  0$  satisfy the inequality
$$
{d\lambda (t)\over dt} \le  - b\lambda (t) +de^{-nt}\ ,\quad t \ge
t_{0}\ ,\quad b > 0\ ,\quad d \ge  0\ ,\quad \lambda (t_{0}) =
\lambda _{0}\ .
$$
 Then $\lambda (t) \rightarrow  0$  as $t \rightarrow  \infty $
 and, also:
\item{$(1)$}  If $n \le  ab$ and $\lambda _{0} \le
c_{0}db^{-1}e^{-nt_{0}}$,  then, for all $t \ge  t_{0}$,
$$
\lambda (t) \le  {c_{0}d\over b} e^{-nt}\ .
$$
\item{$(2)$}
If $n \le  ab$  and $\lambda _{0} > c_{0}db^{-1}e^{-nt_{0}}$,
  then,  for all $t \ge  t_{0}$,
$$
\eqalign{&\lambda (t) \le  \lambda _{0}e^{-ab(t-t_{0})}\ ,\quad
 t_{0} \le  t \le \ot\cr
&\lambda (t) \le  {c_{0}d\over b} e^{-nt}\ ,\quad t >\ot\ ,\quad
\ot= {\ln c_{0}d(b\lambda _{0})^{-1}-abt_{0}\over n-ab} \ .\cr}
$$
\item{$(3)$} If $n > ab$,  then, for all $t \ge  t_{0}$
$$
\lambda (t) \le  Ce^{-ab(t-t_{0})}\ ,\quad
C = \max \{\lambda _{0},c_{0}db^{-1}e^{-nt_{0}}\}\ .
$$

The asymptotic estimate is
$$
\lambda (t) = O(e^{-\beta t})\ ,\quad \beta  = \min \{ab,n\}
$$
as $t \rightarrow  \infty$.
\bigskip

\sect{\bf 2. The Evolution Equations with Uniformly Monotone Operators}

Let $B$ be a real uniformly  convex  and  uniformly  smooth  (reflexive)
Banach space, $(\varphi ,y)$ the pairing (dual product) between
elements $\varphi  \in  B^*$  and
$y \in  B$, $\|\cdot\|$ and $|\cdot| $ the norms
in $B$ and $B^*$, respectively.

We consider stabilization and estimates for the stabilization rate
of
evolution equations of type  (0.1)  where  $A(t)$ is  a  nonlinear
monotone
operator from $B$ to $B^*$ for all $t \ge  t_{0}$.
Let  us  suppose  first  that  the
operator $A(t)$ tends to its limit $A$, the function $f(t)$ tends
to $f$ as $t \rightarrow  \infty $,
limit equation
$$
Ax = f\ ,\quad A: B \rightarrow  B^*\ ,\quad f \in  B^*\ ,\quad x \in  B
\eqno(2.1)$$
 has the solution $\ox$  and scalar functions $\omega _{1}(t)$ and
$\omega _{2}(t)$ are such that
$$
\big|A(t)x-Ax\big|\le  \omega _{1}(t)g\big(\|x\|\big) \ ,\quad\big\|
f(t)-f\big\|  \le  \omega _{2}(t)\ ,\eqno (2.2)$$
 where $\omega _{1}(t) \rightarrow  0$, $\omega _{2}(t) \rightarrow
 0,$  and $g(\xi )$  is  a  continuous  and  nonnegative
function for $\xi  \ge  0$ and bounded (on the bounded sets).

We should remark that relation (0.2), in fact,  reduces  the  equation
(0.1) to a single space $(B^*$ or $H)$ and it removes the principal
difficulties
connected with cardinal distinctions of $B$ and $B^*$ spaces [17,22].

The first difficulty in our situation is that we can no longer use the
traditional Lyapunov functional $V_{1}(x) = 2^{-1}\|x-\ox \|^{2}$.
 Indeed, if we  want  to
decrease functional $V_{1}\big(x(t)\big)$ for each value $t \ge
 t_{0}$, we should suggest
$$
\big(\grad V_{1}\big(x(t)\big),A(t)x(t)\big) \ge  0\ ,
\eqno(2.3)$$
 because in Banach space (see, for instance [2])
$$
{dV_{1}\big(x(t)\big)\over dt} = \left(U\big(x(t)-\ox\big),{dx(t)\over
dt}\right)
$$
 and
$$\grad V_{1}(x) = U(x-\ox )\ .
$$
 Here $U: B \rightarrow  B^*$ is a normalized duality mapping.
In the  case  of $B \neq  H$,
(2.3) makes no sense since $A(t)x(t) \in  B^*$ as well.

The second difficulty is a necessity to produce  a  suitable  equation
(0.1) so far as in (0.1) $x^\prime (t) \in   B$.
These  above  two  difficulties  are
closely connected.

Consider the evolution equation
$$
{d\varphi (t)\over dt} + A(t)x(t) = f(t)\ ,\quad x(t_{0}) = x^{0}\
,\quad t_{0} < t < \infty \ ,\quad \varphi (t)= Ux(t)
\eqno(2.4)$$
 in the dual space $B^*$.
Let us assume operator $A(t)$ is  uniformly  monotone,
i.e., for all $t \ge  t_{0}$, $x \in  B$, $y \in  B$,
$$\big(A(t)x-A(t)y,x-y\big) \ge  c(t)\psi _{1}\big(\| x-y\|\big)\ ,
\quad c(t) \ge  b > 0 \ ,\eqno (2.5)$$
where $\psi _{1}(\xi )$ is a continuous positive function for
all $\xi  \ge  0$, $\psi _{1}(0) = 0,$ and
$$
\overline{\lim_{\xi\to\infty}}  {\psi _{1}(\xi )\over \xi } = \infty\  .
\eqno(2.6)$$
Let us introduce the Lyapunov functional in Banach space by formula
$$
V(\varphi ,y) = 2^{-1}\big(|\varphi | -2(\varphi ,y)+\| y\| \big)
 \ ,\quad\varphi  = Ux\ .\eqno           (2.7)
$$
 It is a nontrivial functional because it is determined  on  the
 elements  of
primary and dual spaces at the same time.
$V(Ux,y)$ is convex,  nonnegative,
differentiable with respect to $\varphi $ and $y$ functionals, and
$$
\eqalignno{&\grad_{\varphi }V(Ux,y) = x-y \in  B\ ,\hbox{
as $y$ is fixed}\ ,\cr
&\grad_{y}V(Ux,y) = Uy-Ux \in  B^*\ , \hbox{as $x$ is fixed\
,}\cr
&V(Ux,y) \ge  2^{-1}\big(\| x\| -\| y\| \big)^{2}& (2.8)\cr
&V(Ux,y) \le  2^{-1}\big(\| x\|+\| y\| \big)^{2}\ .& (2.9)\cr}$$
 Therefore, $\grad_{\varphi }V(Uy,y) = 0$,  $\grad_{\varphi
}V(Ux,x) = 0$, $V(Uy,y) = 0$, $V(Ux,x)=0$,
$\grad_{y}V(Uy,y) = 0,$ $\grad_{y}V(Uy,y) = 0$.
 From (2.8) and (2.9) it  follows  that
$V(Ux,y) \rightarrow  0$ as $\left\Vert{x}\right\Vert \rightarrow
 \infty $ or $\left\Vert{y}\right\Vert \rightarrow  \infty $ and
vice versa.  In Hilbert space $V(Ux,y)
= V_{1}(x,y) = 2^{-1}\left\Vert{x-y}\right\Vert^{2}$.
It turns out that in a Banach  space  the  functional
$V(Ux,y)$  is  connected   with $V_{1}(x,y)$   by   means
  of   its   geometric
characteristics.

Let  us  show  this.  In  [6-8]  we  established   the   following
inequalities for dual mapping:  for all $x,y \in  B$
$$
\eqalignno{&(Ux-Uy,x-y) \le  8\| x-y\|^{2} +
C_{1}\rho _{B}\big(\|x-y\|\big)&(2.10)\cr
&(Ux-Uy,x-y) \le  8|Ux-Uy|^{2} + C_{1}\rho _{B^*}\big(
|Ux-Uy|\big)&(2.11)\cr
&(Ux-Uy,x-y) \ge  (2L)^{-1}\delta _{B}\big(\|x-y\|/C_{2}\big)\ ,
\quad 1 < L < 3.18\ .&(2.12)\cr
&(Ux-Uy,x-y) \ge  (2L)^{-1}\delta _{B^*}\big(|Ux-Uy|/C_{2}\big)
&(2.13)\cr
&|Ux-Uy|\le  C_{2}g^{-1}_{B^*}\big(2C_{2}L
\|x-y\|\big)\ ,\quad g_{B}(\epsilon ) = \delta _{B}(\epsilon )
/\epsilon
&(2.14)\cr
&\|x-y\|\le  C_{2}g^{-1}_{B}\big(2C_{2}L|Ux-Uy|\big)\ ,\quad
g^{-1}_{B}g_{B}(\epsilon ) = \epsilon  .
&(2.15)\cr}$$
 In (2.10)-(2.15) $\rho (\tau)$  is  the modulus  of  smoothness,
 $\delta (\epsilon )$  is  the modulus  of
convexity of space $B$, the constants $C_{1}$ and $C_{2}$, in general,
  depend  on $\left\Vert{x}\right\Vert$
and $\left\Vert{y}\right\Vert$.  However, if $\left\Vert{x}\right\Vert
 \le  R, \left\Vert{y}\right\Vert \le  R$, then

$$
C_{1} = 8\max \{L,R\},\quad C_{2} = 2\max \{1,R\},\quad 1 < L < 3.18.
$$
 In this case (2.10) and (2.12) are the  quantitative  description of a
 well
known mathematical fact:  a duality mapping is uniformly continuous on
each
bounded set in a uniformly smooth Banach space and uniformly monotone
 on each
bounded set in a uniformly convex Banach space.  The estimates
 (2.10)-(2.13)
are derived from the parallogram inequalities in Banach space [8]
$$
2\|x\|^{2} + 2\|y\|^{2} - \|x+y\|^{2} \le  4\|x-y\|^{2} + C_{1}
\rho _{B}(\|x-y\|)$$
$$
2\|x\|^{2} + 2\|y\|^{2} - \|x+y\|^{2} \ge  (4L)^{-1}\delta _{B}
(\|x-y\|/C_{2})$$
 and (2.14), (2.15) are obtained from (2.12), (2.13) (cf.\ [21]).

Let $y$ be an arbitrary fixed point in B.   By  virtue  of  the
 uniform
convexity of $V_{1}(x,y)$ we have
$$
\|2^{-1}(x+y)\|^{2}-\|x\|^{2} \ge  (Ux,y)-(Ux,x) + (2L)^{-1}
\delta _{B}\big(\|x-y\|/2C_{2}\big).
$$
 Now using the identity $(Ux,x) = \left\Vert{x}\right\Vert^{2}$ the
 following relation
$$
(Ux,y) \le \|2^{-1}(x+y)\|^{2}-(2L)^{-1}\delta _{B}\big(\|x-y\|/2C_{2}
\big)
$$
 is obtained.  So far as $g_{B}(\epsilon )$ is a nondecreasing function,
 it is $\delta _{B}(\epsilon /2) \le
2^{-1}\delta _{B}(\epsilon )$.  Hence, one can write
$$
\eqalign {V(\varphi ,y) &= 2^{-1}\|x\|^{2} - (Ux,y) + 2^{-1}\|y\|^{2}
\ge\cr
&\ge  2^{-1}\|x\|^{2} + 2^{-1}\|y\|^{2} -\|2^{-1}(x+y)\|^{2}+
(2L)^{-1}\delta _{B}\big(\|x-y\|\big)/C_{2} \ge\cr
&\ge  (4L)^{-1}\delta _{B}\big(\|x-y\|/C_{2}) + (2L)^{-1}\delta _{B}
\big(\|x-y\|/2C_{2}\big) \ge\cr
&\ge  L^{-1}\delta _{B}\big(\|x-y\|/2C_{2})\cr}.$$
 On the other hand, by the convexity of $V(\varphi ,y)$
$$
V(\varphi ,y) \le  V(Uy,y) + (Ux-Uy,x-y) = (Ux-Uy,x-y)
\eqno(2.16)$$
 holds.  Finally, that implies
$$
V(Ux,y) \ge  L^{-1}\delta _{B}\big(\|x-y\|/2C_{2}\big)
\eqno(2.17)$$
 and
$$
V(Ux,y) \le  8\|x-y\|^{2} + C_{1}\rho _{B}\big(\|x-y\|\big) =
 \mu _{B}\big(\|x-y\|\big) .
\eqno(2.18)$$
 $\left\Vert{x-y}\right\Vert$ and $B$ in the last two inequalities
 can be replaced by $|Ux-Uy|$  and $B^*$.

The functions $\delta _{B}(\epsilon )$  and $\rho _{B}(\tau)$  are
 strongly  increasing,  therefore
yielding
$$
\eqalignno { &\|x-y\|  \le  2C_{2}\delta ^{-1}_{B}\big(LV(Ux,y)\big)
&(2.19)\cr
&\|x-y\|\ge  \mu ^{-1}_{B}\big(V(Ux,y)\big) .
&(2.20)\cr}$$

There  is  a  connection  between  the  functional $V(Ux,y)$  and
   the Young-Fenchel transformation.   Indeed,  let $\phi (x)$  be  a
  function  in B.
Consider the conjugate functional $\phi^*(x^*)$ in space $B^*$
 determined by
$$
\phi^*(x^*) = \sup_x\big((x^*,x)-\phi (x)\big) .
$$
 If $\phi (x)$ is differentiable and $x^* = \phi'(x)$, then
$$
\phi^*(x^*) = (x^*,x)-\phi (x)
$$
 or, otherwise, if $x^* \neq  \phi'(x)$, then
$$
\phi^*(x^*) > (x^*,x)-\phi (x).
$$
 Consequently, the functional
$$
V(x^*,y) = \phi^*(x^*)-(x^*,y)+\phi (y) \ge  0
$$
 for all $x^* \in  B^*$, $y \in  B$ and $V(x^*,y) = 0$ only if
 $x^* = \phi'(y)$. $\phi^*(x^*)$ is  well
known to be a convex.  Let $\phi^*(x^*)$ be a differentiable.
Then  we  have  for
each fixed $y \in  B$,
$$
V'_{x^*}(x^*,y) = \phi^{*\prime}(x^*) - y \in  B\ .
$$
 This shows the functional $V(x^*,y)$ is a convex and its point
of  minimum  is
determined by the relation
$$
\phi^{*\prime}(x^*) = y \ ,
\eqno(2.21)$$
 i.e., $x^* = [\phi^{*\prime}]^{-1}y$ (if $\phi^{*\prime\prime}(x^*)
\neq  0)$.  By analogy, one can write
$$
V'_y(x^*,y) = \phi ^\prime (y) - x^* \in  B^*
$$
 for each fixed $x^* \in  B$.  It gives
$$
\phi ^\prime (y) = x^*\ ,
\eqno(2.22)$$
 i.e., $y = [\phi ^\prime ]^{-1}x^*$ (if $\phi ^{\prime\prime}(y) \neq
 0)$.  So, at the minimum point $(\ox^*,\oy )$
$$\phi ^\prime \phi^{*\prime}(\ox^ *) = \ox^*\ ,\quad \phi^{*\prime}
\phi (\oy ) = \oy
$$
 hold.

The problem of global minimization of $V(x^*,y)$ in the form
$$\oV (x^*) \rightarrow\min_{x^*\in B^*}$$
 is formulated in [13], where
$$\oV(x^*) = \min_{y\in B}V(x^*,y)\ .
$$
 If $\phi^*(x^*)$ is differentiable then
$$
\oV^\prime _{x^*}(x^*) = \phi^{*\prime}(x^*)-[\phi ^\prime ]^{-1}(x^*)
$$
 and the condition of minimum of $V(x^*)$ is written as
$$
\varphi^{*\prime}(x^*) = [\phi ^\prime ]^{-1}(x^*)\ .
\eqno(2.23)$$
 It is easy to see that (2.21), (2.22) and (2.23) coincide.

In terms of the above, the functional $V(x^*,y)$ can be used
as a general
Lyapunov functional in Banach spaces.
Our functional (2.7) is a particular
case of $V(x^*,y)$, for which
$$
\phi (y) = 2^{-1}\|y\|^{2}\ ,\quad \phi^*(x^*) = 2^{-1}|Ux|^{2}\ ,
\quad V'_{x^*}(x^*,y) = x-y\ ,
$$
$$
V'_y(x^*,y) = Uy-Ux\ , \quad\big(\partial _{y}V(\varphi ,y),\partial
_{\varphi }V(\varphi ,y)\big) = -(x-y,Ux-Uy) \le  0
$$
 are satisfied.

Assume now the function $x(t),y(t): [t_{0},\infty ) \rightarrow  B$.
 Consider the  Lyapunov
functional
$$
V\big(\varphi (t),y(t)) = 2^{-1}\big(|\varphi (t)|^{2}-2( \varphi
(t),y(t)) +\|y(t)\|^2\big)\ ,\quad \varphi (t) = Ux(t)\ .
\eqno(2.24)$$

\proclaim Lemma 11.
 Let $x(t)$  and $y(t)$  be strongly  continuous  functions, $y(t)$
 and $\phi (t)$  be weakly differentiable on interval
$[t_{0},\infty )$.   Then, the functional
$(2.24)$ is differentiable and equality
$$
{dV\big(\varphi (t),y(t)\big)\over dt} = \left({d\varphi (t)\over dt},
x(t)-y(t)\right) + \left(Uy(t)-Ux(t),{dy(t)\over dt}\right)
\eqno(2.25)$$
 holds.

 \subheading{Proof {\rm (cf.\ [12])}}
As far as the spaces $B$ and $B^*$ are uniformly smooth,
their norms are differentiable and
$${\partial V\big(\varphi (t),y(t)\big)\over \partial \varphi } =
x(t)-y(t)\ ,\quad {\partial V\big(\varphi (t),y(t)\big)\over
\partial y} = Uy(t)-Ux(t)\ .
$$
 Because of the convexity of $V(Ux,y)$ one has for $t_{0} < s < t$
that
$$\eqalignno{&
V\big(\varphi (t),y(t)\big)\ge  V\big(\varphi (s),y(t)\big) +
\left({\partial V\big(\phi (s),y(t)\big)\over \partial \varphi },
\phi (t)-\phi (s)\right)\cr
\noalign{\hbox{and}}
&V\big(\varphi (s),y(t)\big) \ge  V\big(\varphi (s),y(s)\big) +
\left({\partial V\big(\phi (s),y(s)\big)\over \partial y},y(t)-y(s)
\right)\cr}
$$
 are true.  Therefore, inequalities
$$\eqalign{
V\big(\varphi (t),y(t)\big) \ge  V\big(\varphi (s),y(s)\big)&+
\big(x(s)-y(t),\varphi (t)-\varphi (s)\big)\cr
&
+ \big(Uy(s)-Ux(s),y(t)-y(s)\big)\cr}$$
 and
$$\eqalign{
V\big(\varphi (s),y(s)\big) \ge  V\big(\varphi (t),y(t)\big)&+
\big(x(t)-y(s),\varphi (s)-\varphi (t)\big)\cr
&+ \big(Uy(t)-Ux(t),y(s)-y(t)\big)\cr}
$$
 take place.  Unifying them, one obtains
$$\eqalignno{&
\left(x(t)-y(s),{\varphi (t)-\varphi (s)\over t-s}\right) +
\left(Uy(t)-Ux(t),{y(t)-y(s)\over t-s}\right) \ge\cr
&\qquad \ge  {V\big(\varphi (t),y(t)\big)-V\big(\varphi (s),y(s)\big)\over
t-s}\ge&(2.26)\cr
&\qquad\ge  \left(x(s)-y(t),{\varphi (t)-\varphi (s)\over t-s}\right)
+ \left(Uy(s)-Ux(s),{y(t)-y(s)\over t-s}\right)\ .\cr}
$$
 Since $B$ is a uniformly smooth space, then $U$ is a continuous
mapping.  Using
the conditions of the lemma let us go to the limit in (2.26) as
$s \rightarrow  t$.   We
immediately get
$$\matrix{\displaystyle \left(x(t)-y(t),{d\varphi (t)\over dt}\right)
 + \left(Uy(t)-Ux(t),{dy(t)\over dt}\right) \ge  {dV(Ux(t),y(t))\over
dt} \ge\cr
\displaystyle\ge  \left(x(t)-y(t),{d\varphi (t)\over dt}\right) +
\left(Uy(t)-Ux(t),{dy(t)\over dt}\right)\  .\cr}
$$
 From this we conclude (2.25), i.e.,
$$
{dV(\varphi ,y)\over dt} = \left({\partial V\over \partial \varphi
},{d\varphi \over dt}\right) + \left({\partial V\over
\partial y},{dy\over dt}\right)\ .$$
 The lemma is proven.

Returning to the evolution differential equation (2.4), we will assume
the existence of its weak solution.

\subheading{Definition}  A function $x(t)$ from $[t_{0},\infty )$ to $B$ is
said
  to  be  a  weak
solution of the equation (2.4)  if $x(t)$  is  strongly  continuous  on  set
$[t_{0},\infty ), \varphi (t)$ is weak once-differentiable and for each $f(t)$
 given in $B^*$  and
$x(t_{0}) = x^{0}$ given in $B$ the equality
$$
\left({dUx(t)\over dt},w)\right) +  \big(A(t)x(t),w\big) =  \big(f(t),w
\big)
$$
 is satisfied for all $w \in  B$.

Let us set in (2.24) $y(t) \equiv  \ox $  (cf.\ [11,12].

\proclaim Theorem 1.
Assume the conditions $(2.2), (2.5), (2.6)$ are carried out.
If the weak solution $x(t)$  of $(2.4)$ exists, then the relation
$$\big\|x(t)-\ox \big\|  \rightarrow  0\eqno (2.27)$$
 yields as $t \rightarrow  \infty $.

\pr
 In fact, from (2.4) and (2.25) we have
$$
{dV\big(\varphi (t),\ox\big )\over dt} = \left({dUx(t)\over dt},x(t)
-\ox\right) = - \big(A(t)x(t)-A(t)\ox,x(t)-\ox\big)$$
$$
+ \big(A\ox-A(t)\ox,x(t)-\ox\big)-\big(A\ox-f(t),x(t)-\ox\big)\  .
$$
 Utilizing the condition of uniform monotonicity (2.5) one can obtain
$${dV\big(\varphi (t),\ox\big)\over dt} \le  -c(t)\psi _{1}
\big(\|x(t)-\ox \|\big)   + \big(\omega _{1}(t)g(\|\ox\|)+
\omega _{2}(t)\big)\|x(t)-\ox \|\ .
$$
 Because of (2.6)  the  functional $V\big(\varphi (t),\ox\big)$
is  bounded  (this  fact  can
be proved easily from the contrary).
Let $V\big(\varphi (t),\ox\big) \le  R_{0}$.
Then $\big\Vert x(t)\big\Vert \le  \|  \ox \|
+ 2\sqrt{R_0}= R$ follows from (2.8).  Now we can use
(2.20) with $C_{1} = \max  \{L,R\}$.  We get
$$
{dV\big(\varphi (t),x\big)\over dt} \le  -c(t)\psi
\big(V(\varphi (t),\ox)\big) + \gamma (t)\ ,\qquad t_{0} \le
t < \infty\ ,\quad V\big(\varphi (t_{0}),\ox\big) = V_{0}
$$
 where
$$\gamma (t) = \big(\omega _{1}(t)g(\|\ox\|)+\omega _{2}(t)\big)
(R+\|\ox
\|) \eqno             (2.28)
$$
$$
\psi (V) = \psi _{1}\big(\mu ^{-1}_{B}(V)\big)\ ,\quad
C_{1} = \max \{L,R\}\ .
\eqno(2.29)$$
 Denoting $\lambda (t) = V\big(Ux(t),\ox\big)$ we come to the
differential inequality
$$
{d\lambda (t)\over dt} \le  -c(t)\psi (\lambda ) + \gamma (t)\ ,
\quad\lambda (t_{0}) = V_{0}\ ,\quad t \ge  t_{0}\ .
\eqno(2.30)$$
 By Lemma 1 we assert that $V\big(Ux(t),\ox\big)\rightarrow  0$
as $t \rightarrow  \infty $  because $\gamma (t) \rightarrow  0$  and
$c(t)\ge  b$.
Finally, (2.27) follows from (2.19).  The theorem is proven.

\proclaim Theorem 2.   Assume $(2.2), (2.5)$ and $(2.6)$ exist.
Then  the  solution
$x(t)$  of differential equation $(2.4)$ is bounded
$\big(\Vert x(t) \Vert \le  R\big)$,  and
\item{\rm(a)} in the conditions of Lemma $2, x(t)$  tends to  $\ox$
 and  is  estimated  as
$$\big\|x(t)-\ox \|\le  2C_{2}\delta ^{-1}_{B}\big(Lu(t,C)\big)$$
for all $t \ge  t_{0}$;
\item{\rm(b)} in the conditions of Lemma 3, $x(t)$  tends to $\ox$
   and either
$$\eqalignno{&
\big\|x(t)-\ox\|\le  2C_{2}\delta ^{-1}_{B}\left(L\left(\psi ^{-1}
\left( c_{0}{\gamma (t)\over c(t)}\right)\right)\right)\ ,
\qquad t \ge  t_{0}\cr
\noalign{\hbox{or}}
&\big\|x(t)-\ox\|\le  2C_{2}\delta ^{-1}_{B}\big(Lu(t,C)\big)\ ,
\quad t_0 \le  t<\ot\cr
&\big\|x(t)-\ox\|\le  2C_{2}\delta ^{-1}_{B}\left(L\left(\psi ^{-1}
\left( c_{0}{\gamma (t)\over c(t)}\right)\right)\right)\ ,
\qquad t \ge \ot\cr}$$
 holds.

 {\sl Here
$$
u(t,C) = \phi ^{-1}\bigg[ \phi (C) - a \int^{t}_{t_{0}}c(\tau)
d\tau\bigg]
$$
$$
C = \max  \left\{V_{0},\psi ^{-1}\left( c_{0} {\gamma (t_{0})\over
c(t_{0})} \right)\right\}\ ,\qquad C_{2} = \max \{L,R\}\ ,
$$
 functions $\gamma (t)$  and $\psi (V)$  is determined by
$(2.28),(2.29)$, $\ot$   is a unique  root
of scalar equation}
$$
\phi ^{-1}\bigg[\phi (C) - a \int^{t}_{t_{0}}c(\tau)d\tau \bigg]
= \psi ^{-1}\bigg(c_{0}{\gamma (t)\over c(t)}\bigg)\ .
$$

The proof is obtained from the inequality (2.30), Lemmas 2 and  3,
and
estimate (2.19).
So far as $c(t) \ge  b$ and
$\psi (\lambda )/\lambda  \rightarrow  \infty $ as
$\lambda  \rightarrow  \infty $, the estimates
of convergence rate are determined also by Lemmas 5 and 7 with $m = 0$.

\proclaim Theorem 3.
Assume conditions $(2.2), (2.5)$ and $(2.6)$ are satisfied.
  Then  the
solution $x(t)$  of the differential equation
$$
{dUx(t)\over dt} + \alpha (t)\big(A(t)x(t)-f(t)\big) = 0\ ,
\qquad x(t_{0}) = x^{0}\ ,
\qquad t_{0} \le  t < \infty \ ,
\eqno(2.31)$$
 where
$$
\alpha (t) \rightarrow  0\ ,\qquad \int^{t}_{t_{0}}\alpha
(\tau)d\tau \rightarrow  \infty \ ,\quad{\rm  as }\quad t
\rightarrow  \infty  \ ,
$$
 is bounded $\big(\Vert{x(t)}\Vert \le  R\big)$  and the
assertions of Theorem 2 are proven with
$$
u(t,C) = \phi ^{-1}\bigg[\phi (C) - a \int^{t}_{t_{0}}
c(\tau)\alpha (\tau)d\tau \bigg]
$$
$$
C = \max  \left\{V_{0},\phi ^{-1}\left( c_{0} {\gamma (t_{0})\over
c(t_{0})} \right)\right\}\ ,\qquad C_{2} = \max \{L,R\} \ .
$$

Consider equation (2.1) with a uniformly monotone operator
$$
Ax = Fx + \alpha Ux\ ,\qquad \alpha  =\const > 0\ ,
$$
 where $F$ is a proper monotone operator, i.e., it satisfies
$$
(Fx-Fy,x-y) \ge  0
$$
 and there is no  intensification  of  this  condition.
 The  corresponding
evolution differential equation (2.4) is
$$
{dUx(t)\over dt} + Fx(t) + \alpha Ux(t) = 0\ ,\quad
x(t_{0}) = x^{0}\ ,\quad t \ge  t_{0}\ .
\eqno(2.32)$$
 An equation of this kind has been studied previously in Hilbert
spaces  and
in Banach spaces $B$, when operator $A$ is accretive and acts
from $B$ to $B$.  For
the solution $x(t)$ of the equation
$$
{dx(t)\over dt} + Fx(t)+\alpha x(t) = 0\ ,\quad x(t_{0}) =
x^{0}\ ,\quad t \ge  t_{0}
$$
 the estimate
$$\big\|x(t)-\ox \big\|\le  \big\|x(t_{0})-\ox \big\|
e^{-  \alpha (t-t_0)}$$
 was obtained, where $F\ox+ \alpha \ox   = 0$.
We may establish a similar  estimate  for
equation (2.32).

\proclaim Theorem 4.
Suppose that the weak solution $x(t)$  of the equation $(2.32)$
exists, $Fx$  is a proper monotone operator, $\ox$   is the
 solution  of  equation
$Fx + \alpha Ux = 0$.   Then $x(t)$  approaches $\ox$
strongly as $t \rightarrow  \infty $  and
$$
\big\|x(t)-\ox  \big\|\le  2C_{2}\delta ^{-1}_{B}(LV_{0}
e^{-\alpha (t-t_0)})\ .
$$

\pr  Using property (2.16), we have
$$
\eqalign{{dV(Ux(t),x)\over dt}&= -\big(Fx(t)+\alpha (t)Ux(t),
x(t)-\ox\big)\cr
&= -\big(Fx(t)-F\ox  , x(t)-\ox  \big)-\alpha\big (Ux(t)-U\ox
  ,x(t)-\ox\big  )\cr
&\le  - \alpha V\big(Ux(t),\ox\big)\ .\cr}
$$
 From here the inequality
$$
V\big(Ux(t),\ox\big) \le  V(Ux^{0},\ox)e^{-\alpha (t-t_0)}
= V_{0}e^{-\alpha (t-t_0)}
$$
 is valid.  And the final estimate follows from (2.19).
\hfill$\square$

\medskip
Let us set in (2.25) $y(t) \equiv  0.$  We obtain the identity
$$
{dV(\varphi (t))\over dt} = {1\over 2} {d\big\|x(t)\big\|^2\over dt}
= \left({dUx(t)\over dt},x(t)\right)\ .
\eqno(2.33)$$
 It allows us to produce the statement which is analogous  to
[15]  as  the
corollary above.

\proclaim Theorem 5.    (a)  Let  the  assumptions $\big(A(t)x,x\big)
\ge  0$   and $f(t) \in
L_{2}(t_{0},T;B^*)$  be satisfied.  Then $x(t) \in
L_{\infty }(t_{0},\infty ;B)$.\hb
\indent (b) If $\big(A(t)x,x\big) \ge  \psi \big(\|x\|\big )$  and
$f(t) \in  L_{1}(t_{0},T,B^*)$,  then $x(t) \rightarrow  0$   in $B$
 as $t \rightarrow  \infty $.

 From (2.33) we obtain
$$\big\|  x(t)\big\| ^{2}_{B} \le  \big\|x(s)\big\|^{2}_{B}
- 2 \int^{t}_{s}\big(A(\tau)x(\tau)-f(\tau),x(\tau)\big)d\tau\ ,
$$
 which is similar to Hilbert space.  It is obvious that
$$
\eqalignno{&\left({dUx(t)\over dt},x(t)\right) =
\left(Ux(t),{dx\over dt}\right)\cr
\noalign{\hbox{and}}
&\left({d^{2}Ux(t)\over dt^{2}},x(t)\right)   = \left(
 Ux(t), {d^{2}x\over dt^{2}}\right)\ .\cr}
$$
\bigskip

\sect{3. Evolution Equations with Proper Monotone Operators}

Evolution equations
$$
{dx(t)\over dt} + Ax(t) = f\ ,\qquad x(t_{0}) = x^{0}\ ,\qquad
t \ge  t_{0}
$$
with proper monotone operators are unstable with respect
to  perturbations
of the operator $A$ and right-hand part $f$.  It is known
that one cannot prove
a convergence of solutions $x(t)$  to  the  solution  $\ox$
   of  the  stationary
equation (2.1) in the class of nonlinear  operators  even
 for  unperturbed
equations, the more so if the equations are of the kind (2.4).
Therefore,  in
such a situation one has to apply a stabilization operator.
In the case of
$A: B \rightarrow  B$ being an accretive operator in Banach
space (or monotone  operator
in Hilbert space), a one-parameter operator $\alpha (t)I$
has  been  used  for
stabilization where $\alpha (t) \rightarrow  0$ as
$t \rightarrow  \infty $ and  $I$  is  the  identical  operator.
Under these conditions a convergence of solutions of the evolution
equation
$$
{dx(t)\over dt} + A(t)x(t) + \alpha (t)x(t) = f(t)
$$
 has been proven in [2,12,14].

In the present section we shall consider the evolution equation
$$
{dUx(t)\over dt} + A(t)x(t) + \alpha (t)Ux(t) = f(t)\ ,\qquad
 x(t_{0}) = x^{0}\ ,\qquad
 t \ge  t_{0}\ .
\eqno(3.1)$$

Let us introduce an auxiliary equation
$$
Ay(t) + \alpha (t)Uy(t) = f\  .
\eqno(3.2)$$
 Assume the solution $y(t)$ to be continuous and strongly
 differentiable  for
all $t \ge  t_{0}$.

\proclaim Theorem 6.   Assume the following:
\item{\rm (i)}
   A weak solution of $(3.1)$ exists and is bounded,
$\big\Vert{x(t)}\big\Vert \le  R_{1}$;
\item{\rm (ii)}
  Operator $A(t)$  tends to limit one $A$  as $t \rightarrow
\infty $,  and $f(t)$  tends to  limit
function $f$  with the estimates (2.2);
\item{\rm (iii)}
 The solution $\ox$   of stationary equation $Ax = f$  exists;
\item{\rm (iv)}
 $\alpha (t)$  is a positive, continuous and differentiable function
satisfying
conditions such as $t \rightarrow  \infty $:
$$
\alpha (t) \rightarrow  0\ ,\quad \int^{t}_{t_{0}}\alpha (\tau)d\tau
\rightarrow  \infty \ ,\quad {\big\| \alpha ^\prime (t)\big\|
\over \alpha ^{2}(t)} \rightarrow  0\ ,\quad
{\omega _{1}(t)+\omega _{2}(t)\over \alpha (t)} \rightarrow  0\ ;
$$
\item{\rm (v)}
 Let also either $\delta (\epsilon ) = O(\epsilon ^{2})$  and
$y(t)$  be differentiable;  or
\item{\rm (vi)}
$ \delta (\epsilon ) = 0(\epsilon ^{p})$, $p > 2$,  and
$$
\left({dUy(t)\over dt},{dy\over dt}\right) \ge\left\|   {dy\over dt}
\right\|  \psi \left(\left\Vert{dy\over dt}\right\Vert\right)\ ,
$$
 where function $\psi (\xi )$   is  strictly  increasing  and
 continuous  for $\xi >0$,
$\psi (0)=0$,  and operator $Ax$  is  monotone  and  strongly
 differentiable  with
respect to $x$.  Then
$$\lim_{t\to\infty}\big\|x(t)-\ox  \big\|   = 0\ .
$$

\pr  It is known [3] that $y(t) \rightarrow  \ox$   and $\big\|y(t)
\big\|\le  R_{2}$.  On account of
(2.2), (2.12), (3.2), Lemma 11 and the monotonicity of $A(t)$, we have
$$\eqalignno{&
{dV\big(Ux(t),y(t)\big)\over dt} = \left({d\varphi \over dt},
{dV\over d\varphi }\right) + \left({dV\over dy},{dy\over dt}\right)
=& (3.3)\cr
&\eqalign{= - \big(A(t)x(t)&-A(t)y(t),x(t)-y(t)\big) -
\alpha (t)\big(Ux(t)-Uy(t),x(t)-y(t)\big)-\cr
&- \big(Ay(t)+\alpha Uy(t)-f,x(t)-y(t)\big) +
\big(f(t)-f,x(t)-y(t)\big)+\cr
&+ \big(Ay(t)-A(t)y(t),x(t)-y(t)\big) +
\big(Uy(t)-Ux(t),dy(t)/dt\big)\le\cr}\cr
&\eqalign{\le  - \alpha (t)(2L)^{-1}&\delta _{B}\big(\|x(t)-y(t)\|
/C_{2}\big)
 + \big(\omega _{1}(t)g(R_{2})+\omega _{2}(t)\big)\|x(t)-y(t)\|+\cr
&+\big\|Ux(t)-Uy(t)\big\|\big\|dy(t)/dt\big\|\  ,\quad C_{2} =
\max \{1,R_{1},R_{2}\}\ .\cr}\cr}
$$
 Now $\big\|  dy(t)/dt\big\|$   has to be estimated.
It is shown in [5] that
$$
\big\|y(t_{1})-y(t_{2})\big\|\le  R_{3}g^{-1}_{B}\bigg({R_{3}
\big| \alpha (t_{1})-\alpha (t_{2})\big| \over \alpha
(t_{1})}\bigg)\ ,$$
 where $R_{3} = 2\max\big \{1,\| \ox  \|  \big\}$.
   If $\delta _{B}(\epsilon ) = O(\epsilon ^{2})$  and $y(t)$
  is   strongly
differentiable, then
$$
\left\| {dy(t)\over dt}\right\|  \le  R_{4} {\big| \alpha ^\prime
(t)\big| \over \alpha (t)}\  .
$$
 Further, using (2.20) and (2.14) we obtain
$$\eqalign{ {d\over dt}V\big(Ux(t),y(t)\big)&\le
 - \alpha (t)(2L)^{-1}\delta _{B}\big(\mu ^{-1}_B\big(V(Ux(t),y(t))
\big)/C _2\big)+\cr
&\quad+ \big(\omega _{1}(t)g(R_{2})+\omega _{2}(t)\big)(R_{1}+R_{2})
+ \mu _{B}(R_{1}+R_{2})| \alpha ^\prime (t)| /\alpha (t)\ .\cr}
$$
 One can easily see that the function $\psi (\xi ) =
\delta _{B}\mu ^{-1}_{B}(\xi )$  is  continuous  and
positive for all $\xi  > 0$ and $\psi (0) = 0$.
 From the conditions  of  the  theorem
and from Lemma 1,
$$
V\big(Ux(t),y(t)\big) \rightarrow  0\hbox{ as }t \rightarrow  \infty\  .
$$

Further, if condition (vi) is true, then the identity
$$
\left(A'_y{dy\over dt},w\right) + \alpha ^\prime (t)\big(Uy(t),w
\big) + \alpha (t)\left({dUy(t)\over dt},w\right) = 0
$$
 is also true.  Considering this identity under $w = dy/dt$
 and  using  the
monotonicity of operator $A$, one obtains
$$
\alpha (t)\left({dUy(t)\over dt},{dy\over dt}\right) \le  \big|
\alpha ^\prime (t)\big| \big\|y(t)\big\|\left\|    {dy\over dt}\right\|
\  .
$$
 From here and from (vi)
$$
\alpha (t)\psi \left(\left\Vert{dy\over dt}\right\Vert\right)
\le  \big| \alpha ^\prime (t)\big| R_{2}
$$
 follows.  That means
$$
\left\| {dy\over dt}\right\|  \le  \psi ^{-1}\left(R_{2}{\big|
\alpha ^\prime (t)\big| \over \alpha (t)}\right)\ .
$$
 Now the statement $V\big(Ux(t),y(t)\big) \rightarrow  0$
follows again  from  Lemma  1.   If  we
apply estimate (2.19), we conclude $\big\Vert{x(t)-y(t)}\big\Vert
\rightarrow  0$ as $t \rightarrow  \infty $.  According  to
the remark given at the beginning of the proof, we assert
$\lim\limits_{t\to\infty}\big\|x(t)\big\| = 0$. \hfill$ \square $

\subheading{Remark 3}  In [19] the estimate
$$
{1\over 2}{d^{2}\big\|y(t)\big\|^{2}\over dt^{2}} \ge
\left(  {d^{2}y\over dt^{2}},Uy \right) + m \left\|{dy\over dt}
\right\| ^{2}
$$
 is reached under the condition $y^{\prime\prime}(t)
\in  L^{1}(t_{0},\infty ,B)$ and
$$
(Ux-Uy,x-y) \ge  m\|x-y\|^{2}\ ,\qquad m > 0\ .
$$
 This corresponds to $\delta (\epsilon ) = O(\epsilon ^{2})$ if
$\big\Vert{x(t)}\big\Vert \le  R$, $\big\Vert{y(t)}\big\Vert \le  R$.
 From this we
obtain the inequality (cf.\ with condition (vi))
$$
\left({dUy(t)\over dt},{dy\over dt}\right) \ge  m \left\|{dy\over
dt}\right\|^{2}\ .
$$

\subheading{Remark 4}  If in Theorem 5,
$$
\left\| {dy\over dt}  \right\|\bigg/ \alpha (t) \rightarrow
 0\hbox{ as }t \rightarrow  \infty \ ,
$$
 then the assumption  about  the  monotonicity  of  operator
$A(t)$ and  the
differentiability of $A$ could be omitted.
In  this  case  (3.3)  should  be
rewritten as
$$
\eqalign{&\eqalign{{dV\big(Ux(t),y(t)\big)\over dt} \le&
- \big(Ax(t)-Ay(t),x(t)-y(t)\big) - \alpha
(t)\big(Ux(t)-Uy(t),x(t)-y(t)\big)-\cr
&
- \big(Ay(t)+\alpha (t)Uy(t)-f,x(t)-y(t)\big) +
\big(f(t)-f,x(t)-y(t)\big)+\cr
&
+ \big(Ax(t)-A(t)x(t),x(t)-y(t)\big) +
\big(Uy(t)-Ux(t),dy(t)/dt\big)\le\cr}\cr
&\eqalign{
\le  - \alpha (t)(2L)^{-1}\delta _{B}&\big(\|x(t)-y(t)\|/C_{2}\big)
 + \big(\omega _{2}(t)+\omega _{1}(t)g(R_{1})\big)\|x(t)-y(t)\|+\cr
&+ \big\| Uy(t)-Ux(t)\big\|  \big\|  dy(t)/dt\big\|\  .\cr}\cr}
$$

\subheading{Remark 5}  If $\big\| x(t)\big\|\le  R_{1}$ is not
known a priori, the theorem  remains
true but the proof is rather complicated (cf.\ [9]).

\subheading{Remark 6}
The lemmas 2-10 give estimates of a  convergence  rate  of
$\big\Vert x(t)-y(t)\big\Vert$ to zero.

\subheading{Remark 7}
If equation (2.1) has a set of solution $N$, then $x(t)\rightarrow\ox^
 *$,
where $\|\ox^*\|  = \min\limits_{\ox\in N}\|\ox\|$ [3].

\subheading{Remark 8}
Let us determine the  integral  solution  of  (3.4)  as  a
continuous function $x(t) = [t_{0},\infty )$ such that
$$
V\big(Ux(t),y\big) \le  V\big(Ux(s),y\big) + \int^{t}_{s}\big(
f(\tau)-Ay,x(\tau)-y\big)d\tau
$$
 for all $y \in  B$, $t_{0} \le  s \le  t < \infty $ and given
a value of $x(t_{0}) = x^{0}$.   For  this
solution the above theorems are valid [20].

Thus it has been shown that stabilization and the stabilization  rate
of solutions of the evolution equation (3.1),(2.4)(2.31),(2.32) depend  not
only on the structure and smoothness of the problem's operators,  but
 also
on the geometric characteristics of the Banach spaces.

\subheading{Remark 9}
The trajections in dual spaces for problems of  functional
minimization were considered also in [18,9].

\subheading{Remark 10}
As approximate methods, one can use the methods described
in [6,7].
\bigskip
\sect{References}

\item{[1]} Ya.I.~Alber, A continuous regularization of linear
operator  equations
in Hilbert spaces,  Math. Zametki  4 (1968), 503-509 (Russian).
\item{[2]} Ya.I.~Alber, The solution by the  regularization
 method  of  operator
equations of the first kind with accretive operators in Banach  space,
 Differential Equations 11 (1975), 2242-2248.
\item{[3]} Ya.I.~Alber,  The  solution  on  nonlinear  equations
 with  monotone
operators in Banach space,  Siberian Math.\ J.\  16 (1975), 1-8.
\item{[4]} Ya.I.~Alber, Recurrence relations and variational
inequalities,  Soviet
Math.\ Dokl.\ 27 (1983), 511-517.
\item{[5]} Ya.I.~Alber, Iterative  regularization  in  Banach  spaces,
 Izvestiya
Vuzov.\ Math.\ 30 (1986), 1-8.
\item{[6]} Ya.I.~Alber, A.I.~Notik, Geometric properties  of
Banach  spaces  and
approximate methods for solving nonlinear operator  equations,
  Soviet
Math.\ Dokl.\ 29 (1984), 611-615.
\item{[7]} Ya.I.~Alber, A.I.~Notik, On minimization of functionals
 and  solution
of variational inequalities in Banach spaces,  Soviet  Math.\  Dokl.\
 34
(1987), 296-300.
\item{[8]} Ya.I.~Alber, A.I.~Notik, Parallelogram inequalities in
 Banach  spaces
and some properties of  a  duality  mapping,   Ukrainian  Math.\  J.\
40
(1988), 650-652.
\item{[9]} Ya.I.~Alber,  A.I.~Notik,  Iterative  method  for
solving  unstable
variational inequalities on approximately given sets (in preparation).
\item{[10]} Ya.I.~Alber,  S.V.~Shilman,  Recurrent  numerical
and  differential
inequalities, III. \break Preprint, NIRFI, No. 134, 1980.
\item{[11]} H.~Brezis,  On some degenerate nonlinear  parabolic
equations,   Proc.\
Symp.\ Pure Math.\ 18,  Part  I,  Amer.\  Math.\  Soc.,
  Providence,  1970,
28-38.
\item{[12]} F.E.~Browder,  Nonlinear  operators  and   nonlinear
  equations   of
evolutions in Banach space,  Proc.\ Symp.\ Pure Math.\ 18, *Part II,  Amer.
Math.\ Soc., Providence, 1970.
\item{[13]} M.~Dietrich,  On  the  convexification  procedure
for  nonconvex  and
nonsmooth infinite dimensional optimization problems,  J.\  Math.\
 Anal.\
Appl.\ 161 (1991), 28-34.
\item{[14]} M.~Israel, S.~Reich, Asymptotic behaviour of solutions
of a  nonlinear
evolution equation,  J.\ Math.\ Anal.\ Appl.\ 83 (1981), 43-53.
\item{[15]} J.~Kacur, Stabilization of solutions of abstract
parabolic  equations,
 Chechoslovak.\ Math.\ J.\ 30 (1980), 539-555.
\item{[16]} V.~Lakshmikantham,   S.~Leela,     Differential
 and    Integral
Inequalities, Theory and Applications, Academic  Press  New
York  and
London, 1969, Vol.\ I and II.
\item{[17]} J.L.~Lions, Quelques m\'ethodes de r\'esolution
des probl\'emes aux  limites
non lin\'eaires, Dunod Gauthier-Villars, Paris, 1969.
\item{[18]} A.S.~Nemirovskii,  D.B.~Yudin,   Problem  Complexity
 and  Method  in
Optimization, Wiley, Interscience, 1983.
\item{[19]} E.I.~Poffald, S.~Reich, An incomplete Gauchy problem,
 J.\  Math.\  Anal.\
Appl.\ 113 (1986), 514-543.
\item{[20]} E.~Schechter,  Perturbations   of   regularizing
  maximal   monotone
operators,  Israel J.\ of Math.\ 43 (1982), 49-61.
\item{[21]} Zong-Ben Xu, G.F.~Roach,  Characteristic
inequalities  of  uniformly
convex and uniformly smooth Banach spaces,  J.\ Math.\  Anal.\  Appl.\
 157
(1991), 189-210.
\item{[22]} E.~Zeidler,  Nonlinear Functional Analysis and its
 Applications  II/B:
 Nonlinear Monotone Operators, Springer-Verlag, New-York Inc.\ 1990.
\end